\documentclass[letterpaper, 10pt, conference]{ieeeconf}      
\IEEEoverridecommandlockouts                              

\let\labelindent\relax
\usepackage{graphicx}   
\usepackage{epsfig} 
\usepackage{svg}
\usepackage{amssymb}

\usepackage{amssymb} 
\usepackage{caption}
\usepackage{subcaption}
\usepackage{booktabs}
\usepackage{amsmath} 
\usepackage[pdf]{pstricks}
\usepackage{algorithm}
\usepackage[noend]{algpseudocode}
\usepackage{graphicx}
\usepackage{lipsum}%
\usepackage{xcolor}
\usepackage{hyperref}

\usepackage{float}
\usepackage{enumitem}
\usepackage[backend=biber,style=ieee,sorting=none, isbn=false, doi=false, url=false]{biblatex}
\newcommand{\RNum}[1]{\uppercase\expandafter{\romannumeral #1\relax}}
\algdef{SE}[SUBALG]{Indent}{EndIndent}{}{\algorithmicend\ }%
\algtext*{Indent}
\algtext*{EndIndent}
\usepackage[many]{tcolorbox}    	

\tcbset{
    sharp corners,
    colback = white,
    before skip = 0.2cm,    
    after skip = 0.5cm      
}        

\newtcolorbox{boxA}{
    fontupper = \bf,
    boxrule = 1pt,
    colframe = black 
}

\usepackage{color}

\title{\LARGE \bf Connecting the Dots: Context-Driven Motion Planning Using Symbolic Reasoning}
\author{Chris van der Ploeg$^{1,2}$, Michiel Braat$^{1}$, Beatrice Masini$^{3}$, Jochem Brouwer$^{1}$, Jan-Pieter Paardekooper$^{1,4}$ 
\thanks{$^{1}$Netherlands Organisation for Applied Scientific Research, Integrated Vehicle Safety Group, 5700 AT Helmond, The Netherlands.}
\thanks{$^{2}$Eindhoven University of Technology, Dynamics and Control Group, Mechanical Engineering Dept., P.O. Box 513, 5600 MB, Eindhoven, The Netherlands.(e-mail: \href{mailto:c.j.v.d.ploeg@tue.nl}{c.j.v.d.ploeg@tue.nl}).}
\thanks{$^{3}$Netherlands Organisation for Applied Scientific Research, Intelligent Autonomous Systems Group, 2597 AK The Hague, The Netherlands.}
\thanks{$^{4}$Radboud University, Donders Institute for Brain, Cognition and Behaviour, P.O.Box 9010, 6500 GL, Nijmegen, The Netherlands}
}
\date{January 2023}

\begin{document}
\maketitle
\thispagestyle{empty}
\pagestyle{empty}
\begin{abstract}
    The introduction of highly automated vehicles on the public road may improve safety and comfort, although its success will depend on social acceptance. This requires trajectory planning methods that provide safe, proactive, and comfortable trajectories that are risk-averse, take into account predictions of other road users, and comply with traffic rules, social norms, and contextual information. To consider these criteria, in this article, we propose a non-linear model-predictive trajectory generator. The problem space is populated with risk fields. These fields are constructed using a novel application of a knowledge graph, which uses a traffic-oriented ontology to reason about the risk of objects and infrastructural elements, depending on their position, relative velocity, and classification, as well as depending on the implicit context, driven by, e.g., social norms or traffic rules. Through this novel combination, an adaptive trajectory generator is formulated which is validated in simulation through 4 use cases and 309 variations and is shown to comply with the relevant social norms, while taking minimal risk and progressing towards a goal area. 
\end{abstract}
\section{Introduction}
The task of a trajectory generator in an automotive application is to plan safe and comfortable trajectories while ensuring progress towards a goal point~\cite{van_der_ploeg_long_2022}. To achieve this, the vehicle needs to perceive its environment (e.g., objects and the infrastructure), comprehend the current situation, and project into the future. These three elements combined are commonly classified using the definition of situational awareness~\cite{endsley_toward_1995}.  The type of information involved in generating situational awareness can be differentiated into "\textit{explicit}" and "\textit{implicit}" information. By explicit information, we refer to the information which can be extracted directly from sensor measurements, such as the estimated state of the vehicle itself, or the observed state of other objects and the environment. When planning a trajectory, based purely on explicit information, aspects that would require inference to derive usable implicit information are neglected. Examples of these aspects are traffic rules, social norms, or  potentially hidden objects. Although traffic rules and social norms could be noted down explicitly, it often requires situational comprehension to correctly interpret and respect the rules/norms. For implicit information or contextual information, we borrow the definition from~\cite{dey_understanding_2001}, i.e., any information that can be used to characterize the situation of an entity. An entity is a person, place, or object that is considered relevant to the interaction between a user and an application, including the user and applications themselves. In this work, we propose a situational aware trajectory generator that incorporates contextual information, thus making it context-aware.
\subsection{Related work}
Several works use a model-predictive control (MPC) -based trajectory generator, where incorporating objectives explicit observations (road and objects) are either formulated using constraints~\cite{gutjahr_lateral_2017} or through the use of artificial potential fields (APFs) augmented in the objective~\cite{ji_path_2017}. Multiple works collect and make use of contextual information for various purposes. In~\cite{halilaj_knowledge_2021}, a graph-based method is used for identifying context, assessing the difficulty of situations, and trajectory prediction for offline identification of traffic situations. In~\cite{paardekooper_hybrid-ai_2021}, contextual information is employed for evaluating whether an intelligent subsystem of an automated vehicle understands its environment well enough to give reliable outputs. 

In relation to incorporating context in motion planning by, e.g., social norms or traffic rules-based, several works have been proposed in the literature. In~\cite{chen_socially_2017}, a reinforcement learning approach is provided which, during training, employs reward functions inspired by context-driven socially compliant behavior. In \cite{Huang_2019} decision-making for changing lanes is done by a prolog reasoner based on the safety legitimacy and reasonableness of changing lanes. In \cite{Zhao_2017} a knowledge base and a reasoning system are combined to formalize right-of-way rules and make decisions to wait or go in uncontrolled intersections.  In~\cite{xie_distributed_2022} and~\cite{hang_human-like_2021}, the use of risk fields is proposed to induce "human-like" driving and considerations. In these approaches, certain maneuvers (e.g., to make a lane change) are however still largely initiated through discrete decision-making algorithms, which could sacrifice the optimality of the solution. 

In this work, we focus on the comprehension of the current situation using contextual information. Instead of using this comprehension to directly make driving decisions, as is done in \cite{Huang_2019} \cite{Zhao_2017}, contextual information is used to shape the risk field. By shaping this risk field, we do not directly make decisions about what the trajectory should look like, instead, we leave it up to the trajectory generator to optimally weigh risk, progression, and comfort given the constraints.
\subsection{Contributions}
This work contributes to the state-of-the-art in this field through the following aspects:
\begin{itemize}
    \item We formalize the terminology and role of context information for motion planning of vehicles and propose a risk-based knowledge graph that incorporates explicit information and forms the contextual information.
    \item We combine a model-based predictive planner with a knowledge graph. By constructing this hybrid setup, the planner can construct trajectories for complex maneuvers, based on contextual information, driven by, e.g., social norms, traffic rules, or inherent object properties based on their classification and environment.
    \item We evaluate our methodology through four use cases with a total of $309$ variations in simulation, showing the robustness and real-time applicability of our approach.
\end{itemize}
The outline of our work is as follows. First, the problem statement and some required preliminaries are introduced in Sec.~\ref{sec:problemstatement}. Our proposed method for generating situational awareness is introduced and elaborated in Sec.~\ref{sec:method0}. Including this situational awareness into a model-predictive framework is explained in Sec.~\ref{sec:method1}. Subsequently, the proposed approach is demonstrated in Sec.~\ref{sec:simulation} and concluded in Sec.~\ref{sec:conclusion}.
\section{Preliminaries and Problem Statement}\label{sec:problemstatement}
\begin{figure}[t]
    \centering
    \fontsize{9pt}{11pt}\selectfont
    \def\svgwidth{0.4\columnwidth}
\begingroup%
  \makeatletter%
  \providecommand\color[2][]{%
    \errmessage{(Inkscape) Color is used for the text in Inkscape, but the package 'color.sty' is not loaded}%
    \renewcommand\color[2][]{}%
  }%
  \providecommand\transparent[1]{%
    \errmessage{(Inkscape) Transparency is used (non-zero) for the text in Inkscape, but the package 'transparent.sty' is not loaded}%
    \renewcommand\transparent[1]{}%
  }%
  \providecommand\rotatebox[2]{#2}%
  \newcommand*\fsize{\dimexpr\f@size pt\relax}%
  \newcommand*\lineheight[1]{\fontsize{\fsize}{#1\fsize}\selectfont}%
  \ifx\svgwidth\undefined%
    \setlength{\unitlength}{379.84251969bp}%
    \ifx\svgscale\undefined%
      \relax%
    \else%
      \setlength{\unitlength}{\unitlength * \real{\svgscale}}%
    \fi%
  \else%
    \setlength{\unitlength}{\svgwidth}%
  \fi%
  \global\let\svgwidth\undefined%
  \global\let\svgscale\undefined%
  \makeatother%
  \begin{picture}(1,0.89552239)%
    \lineheight{1}%
    \setlength\tabcolsep{0pt}%
    \put(0,0){\includegraphics[width=\unitlength,page=1]{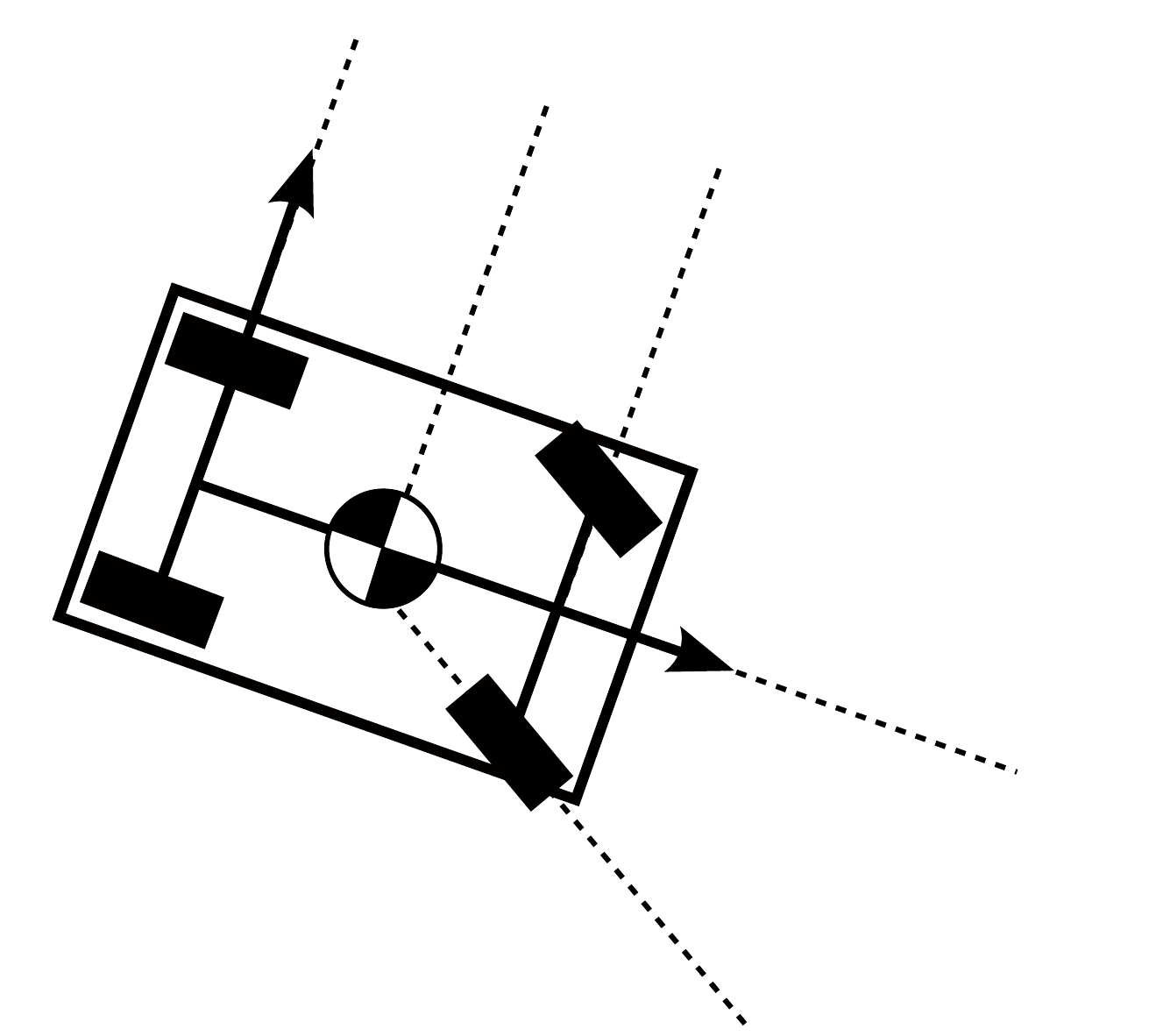}}%
    \put(0.61006759,0.35268583){\color[rgb]{0,0,0}\rotatebox{-1.098766}{\makebox(0,0)[lt]{\lineheight{1.25}\smash{\begin{tabular}[t]{l}$x,v$\end{tabular}}}}}%
    \put(0.17104855,0.67326789){\color[rgb]{0,0,0}\rotatebox{-0.552413}{\makebox(0,0)[lt]{\lineheight{1.25}\smash{\begin{tabular}[t]{l}$y$\end{tabular}}}}}%
    \put(0,0){\includegraphics[width=\unitlength,page=2]{kinematic_model_standalone.pdf}}%
    \put(0.448214,0.82740204){\color[rgb]{0,0,0}\makebox(0,0)[lt]{\lineheight{1.25}\smash{\begin{tabular}[t]{l}$L$\end{tabular}}}}%
    \put(0.3402279,0.78126749){\color[rgb]{0,0,0}\makebox(0,0)[lt]{\lineheight{1.25}\smash{\begin{tabular}[t]{l}$L_r$\end{tabular}}}}%
    \put(0.48387505,0.6924149){\color[rgb]{0,0,0}\makebox(0,0)[lt]{\lineheight{1.25}\smash{\begin{tabular}[t]{l}$L_f$\end{tabular}}}}%
    \put(0,0){\includegraphics[width=\unitlength,page=3]{kinematic_model_standalone.pdf}}%
    \put(0.72958431,0.05240956){\color[rgb]{0,0,0}\makebox(0,0)[lt]{\lineheight{1.25}\smash{\begin{tabular}[t]{l}$\delta$\end{tabular}}}}%
  \end{picture}%
\endgroup%

    \caption{Kinematic bicycle model.}
    \label{fig:3_model}
    \vspace{-0.5cm}
\end{figure}

We borrow the notation and definitions from~\cite{van_der_ploeg_long_2022} to support in step-wise building up to the problem statement of this work. Consider the automated vehicle, which can be described using difference equations, describing the kinematic bicycle model, with discrete sampling time $t_s$, as follows, 
\begin{align}
\begin{aligned}
    x_{k+1} =& x_k+t_sv_k\cos{\theta_k},\quad y_{k+1} = y_k+t_sv_k\sin{\theta_k},\\
    \theta_{k+1} =& \theta_k+t_s\frac{v_k}{L}\tan{\delta_k},\quad v_{k+1} = v_k+t_sa_k,\\
    \delta_{k+1} =& \delta_k+t_s\omega_k.
\end{aligned}\label{eq:vehiclemodel}
\end{align}
where the states $x,y$ denote the planar position of the vehicle, the state $\theta$ denotes the heading, $v$ denotes the longitudinal velocity, $\delta$ denotes the steering angle, $a$ denotes the input longitudinal acceleration, $\omega$ denotes the input steering-rate and $L$ denotes the wheelbase of the vehicle. The vehicle model is also depicted in Fig.~\ref{fig:3_model}~\cite{van_der_ploeg_long_2022}. We would like our trajectory to comply with this dynamical model, as it captures kinematic behavior, suitable to urban situations, and allows constraining the trajectory in states which could benefit safety and comfort (such as the acceleration or the steering rate).  Note that the state is defined from the center of the rear axle of the vehicle.
We hence identify our states and inputs to belong to the compact sets $\mathbf{X},\mathbf{U}$, as follows.
\begin{align}
    \mathbf{X}\!=\!\{\mathbf{x}\in\mathbb{R}^5|\underline{\mathbf{x}}\leq\mathbf{x}\leq\overline{\mathbf{x}}\}, \mathbf{U}\!=\!\{\mathbf{u}\in\mathbb{R}^2|\underline{\mathbf{u}}\leq\mathbf{u}\leq\overline{\mathbf{u}}\}.\label{eq:boundedsets}
\end{align}
The automated vehicle is expected to drive from an initial state given by $(x_0,y_0,\theta_0,v_0,\delta_0)$, of which the planar states $x,y,\theta$ are located on top of the coordinate system (i.e., $x_0,y_0,\theta_0=0$), towards a reference state $\mathbf{r}_N$ which is positioned relative to the initial condition of the vehicle, and can be obtained from a global planner (e.g., Dijkstra or A*, not in the scope of this work). The reference contains a planar position and longitudinal velocity $v$ as a reference, i.e., $(x_N,y_N,\theta_N,v_{N})$.  Here, $N$ denotes the horizon of the trajectory generator to be planned towards the goal state, i.e., the look-ahead time. Throughout this work, we assume that the reference state, $\mathbf{r}_N$, is located in the center of the rightmost lane of the road at a distance of $N\cdot t_s\cdot \bar{v}$, where $t_s$ represents the sampling time between each trajectory element and $\bar{v}$ represents the known maximally allowed velocity in the environment. The vehicle is expected to maneuver in an environment that contains a set of entities, $\mathcal{E}$, of which the position of each entity at the respective time, $k$, i.e., $\mathbf{e}_k\in\mathcal{E}$ is measured with respect to the initial state of the ego-vehicle at the time of planning the trajectory (i.e., $(x_0,y_0,\theta_0)$). The role and hierarchy of entities in this work are further elaborated in Sec.~\ref{sec:method0}. Identically to~\cite{van_der_ploeg_long_2022}, the road markers/boundaries are described through a polynomial description $\mathbf{p}^{(i)}\in\mathcal{E}$, where the superscript denotes the $i$-th polynomial, assumed to be measured with respect to the initial condition of the automated vehicle. The polynomials used in the remainder of this work are, without loss of generality, assumed to be of the third order. The objects, e.g., pedestrians, vehicles, or static obstacles, are defined as $\mathbf{o}^{(j)}_k\forall k\in[0\hdots N]$ with planar state $(x^{o(j)}_k,y^{o(j)}_k,\theta^{o(j)}_k)$, for all observed objects $j$, which are measured and predicted with respect to the initial condition of the vehicle. To include the predicted motion objects in our work, the state of these objects is assumed to be exactly measurable/predictable up to the prediction horizon $N$.

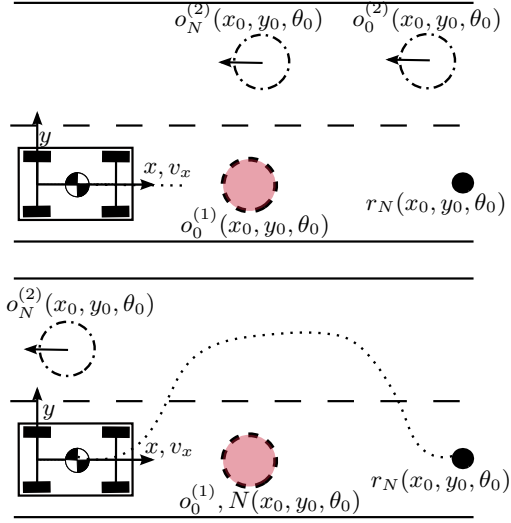
\begin{figure}[t]
    \centering
    \begin{subfigure}{\columnwidth}
    \centering
        \def\svgwidth{0.75\columnwidth}
        \fontsize{9pt}{11pt}\selectfont
        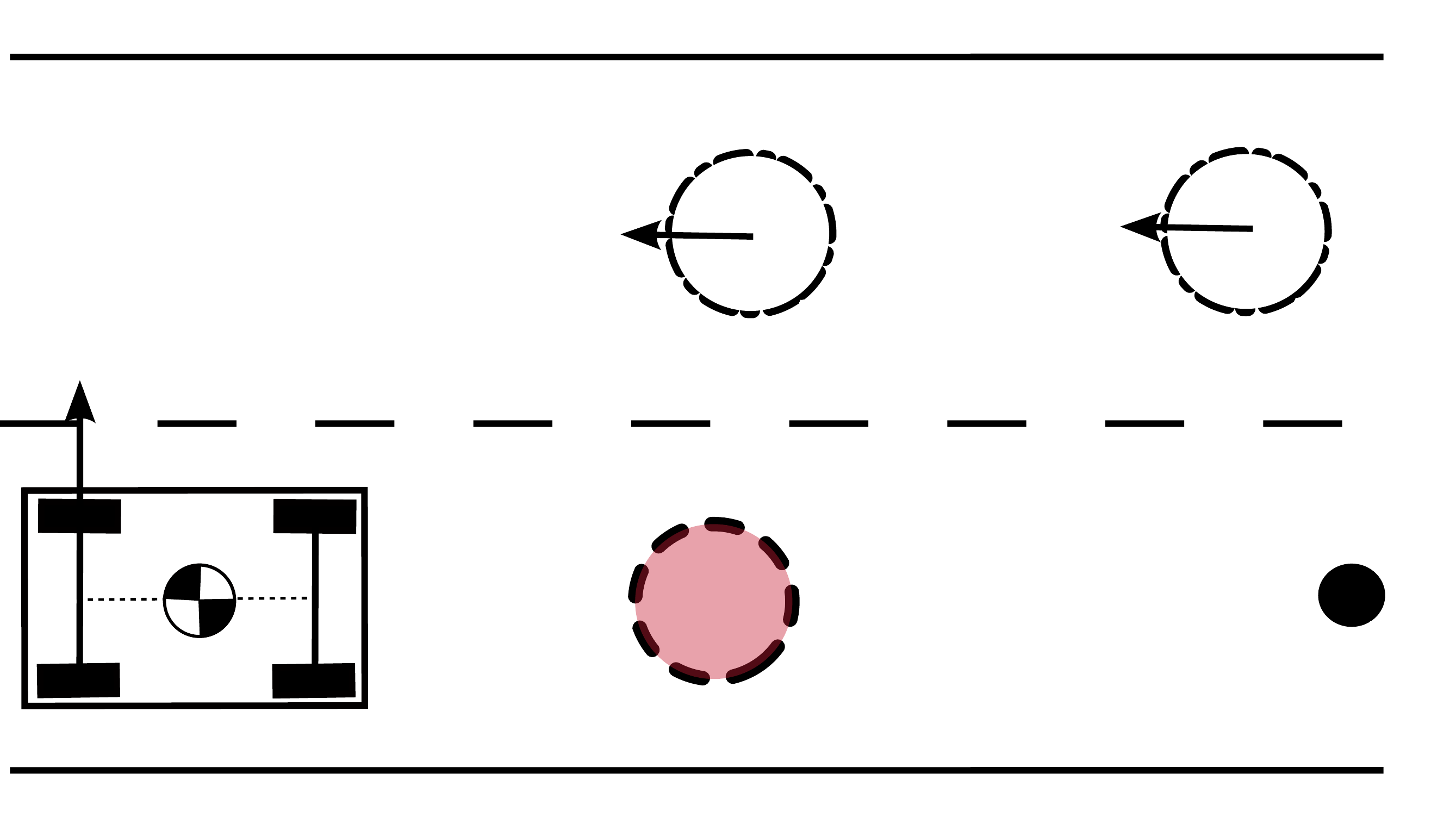
    \end{subfigure}
    \begin{subfigure}{\columnwidth}
    \centering
        \def\svgwidth{0.75\columnwidth}
        \fontsize{9pt}{11pt}\selectfont
\begingroup%
  \makeatletter%
  \providecommand\color[2][]{%
    \errmessage{(Inkscape) Color is used for the text in Inkscape, but the package 'color.sty' is not loaded}%
    \renewcommand\color[2][]{}%
  }%
  \providecommand\transparent[1]{%
    \errmessage{(Inkscape) Transparency is used (non-zero) for the text in Inkscape, but the package 'transparent.sty' is not loaded}%
    \renewcommand\transparent[1]{}%
  }%
  \providecommand\rotatebox[2]{#2}%
  \newcommand*\fsize{\dimexpr\f@size pt\relax}%
  \newcommand*\lineheight[1]{\fontsize{\fsize}{#1\fsize}\selectfont}%
  \ifx\svgwidth\undefined%
    \setlength{\unitlength}{771.02362205bp}%
    \ifx\svgscale\undefined%
      \relax%
    \else%
      \setlength{\unitlength}{\unitlength * \real{\svgscale}}%
    \fi%
  \else%
    \setlength{\unitlength}{\svgwidth}%
  \fi%
  \global\let\svgwidth\undefined%
  \global\let\svgscale\undefined%
  \makeatother%
  \begin{picture}(1,0.56617647)%
    \lineheight{1}%
    \setlength\tabcolsep{0pt}%
    \put(0,0){\includegraphics[width=\unitlength,page=1]{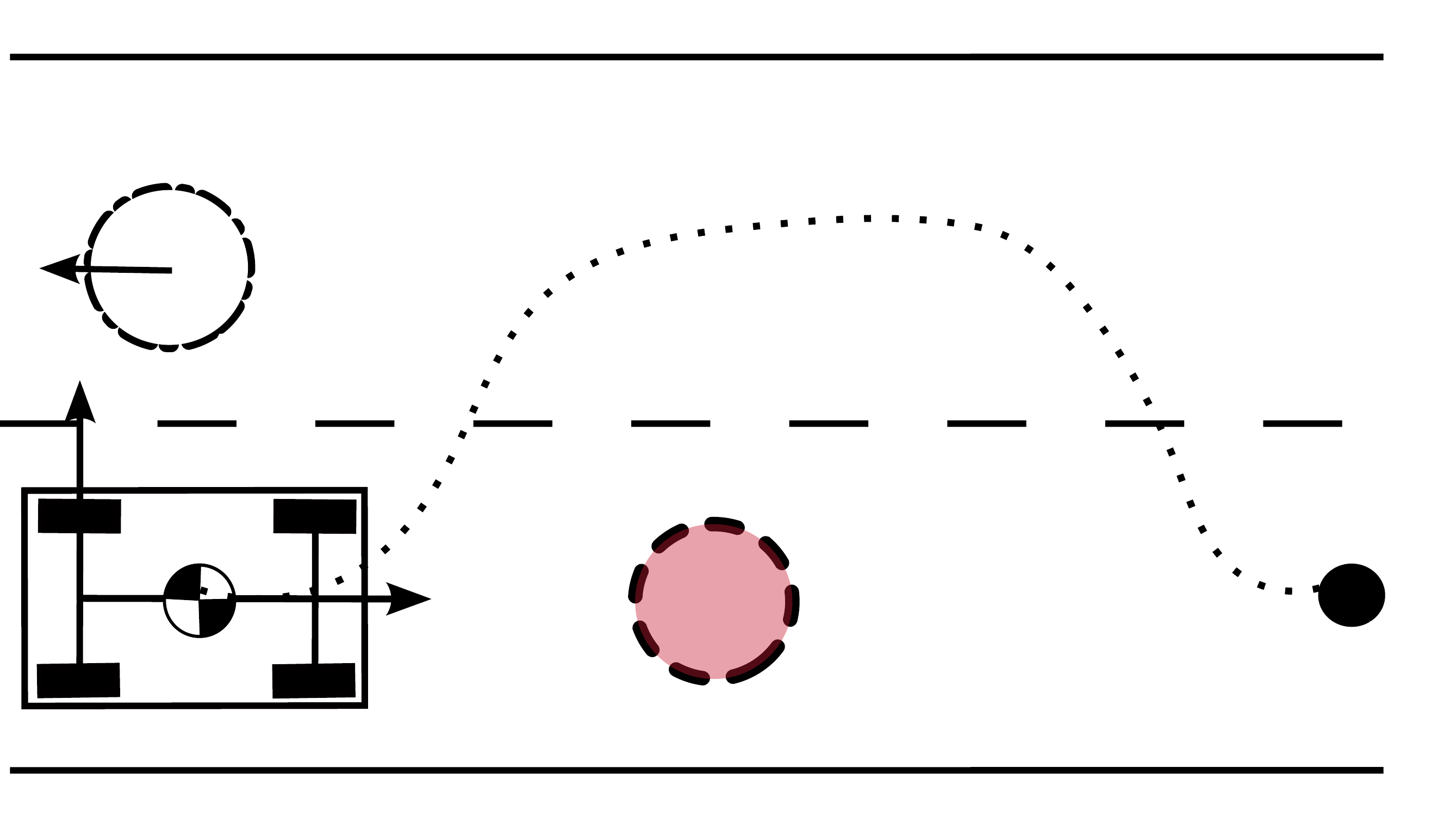}}%
    \put(-0.00921354,0.45979142){\color[rgb]{0,0,0}\makebox(0,0)[lt]{\lineheight{1.25}\smash{\begin{tabular}[t]{l}$o^{(2)}_N(x_0,y_0,\theta_0)$\end{tabular}}}}%
    \put(0.34618751,0.0551346){\color[rgb]{0,0,0}\makebox(0,0)[lt]{\lineheight{1.25}\smash{\begin{tabular}[t]{l}$o^{(1)}_0,N(x_0,y_0,\theta_0)$\end{tabular}}}}%
    \put(0.73940501,0.09793486){\color[rgb]{0,0,0}\makebox(0,0)[lt]{\lineheight{1.25}\smash{\begin{tabular}[t]{l}$r_N(x_0,y_0,\theta_0)$\end{tabular}}}}%
    \put(0.0656545,0.24602961){\color[rgb]{0,0,0}\rotatebox{-0.552413}{\makebox(0,0)[lt]{\lineheight{1.25}\smash{\begin{tabular}[t]{l}$y$\end{tabular}}}}}%
    \put(0.27599637,0.16900119){\color[rgb]{0,0,0}\rotatebox{-1.098766}{\makebox(0,0)[lt]{\lineheight{1.25}\smash{\begin{tabular}[t]{l}$x,v_x$\end{tabular}}}}}%
  \end{picture}%
\endgroup%

    \end{subfigure}
\caption{Two distinct scenarios illustrating the problem statement of our work.}
\label{fig:scenarios}
\vspace{-0.5cm}
\end{figure}

To illustrate the challenges and boil to the main problem statement of this work, we defer the readers' attention to a brief example in Fig.~\ref{fig:scenarios}. Here, the ego vehicle is trying to follow its reference $r_N$ and is approaching a standstill object on its lane, $o_N^{(1)}$. On the opposite lane, $o_N^{(2)}$ is approaching. It is within the line of expectations that the vehicle will stand still behind the static object until $o_N^{(2)}$ has safely passed the ego-vehicle. It is shown in~\cite[Sec. IV.C and IV.D]{van_der_ploeg_long_2022}, that we could fulfill this behavior, given the appropriate tuning of parameters. Now, at the bottom part of Fig.~\ref{fig:scenarios}, the oppositely approaching object has passed, meaning that the ego vehicle would now be able to pass the static object by driving onto the left lane. This infers that the repulsive force of the stationary object is large enough to steer the vehicle around it while crossing the lane toward the left. The 
 work in~\cite[Sec. IV.A, IV.B]{van_der_ploeg_long_2022} shows that the planner would be capable of suggesting such a maneuver. Note, however, that this planner has no understanding that the left lane is intended for traffic driving from the opposite side and, as such, it would be undesirable to remain in that lane once the static object has been passed. In a scenario where $o_N^{(1)}$ is not just a parked object, but a pedestrian crossing the road, it would even be socially undesired or illegal to overtake the pedestrian. Now, consider the case that the static object in Fig.~\ref{fig:scenarios} is not a vehicle or a pedestrian, but a trashcan or a pile of leaves. Clearly, the risk around different classes of an object should induce different shapes of risk fields. The same holds for road marking types, where the risk field of crossing a solid line should be much higher due to its illegality than crossing a dashed line. In summary, the risk field is highly dependent on context and should be time-varying by nature. In summary, a trajectory generator needs to plan a trajectory $\mathcal{T}$ as a function of its states, i.e.,
to drive along the road infrastructure (i.e., follow the lane markings), and interact with other road users (and their predictions) and static elements (e.g., garbage bins, parked vehicles). Moreover, the vehicle is required to abide by traffic rules and comply with social norms. 

\section{Generating Contextual Information}~\label{sec:method0}
Contextual information is formalized using a knowledge graph. The main purpose of the knowledge graph is to aggregate and relate all the information about the current context. Firstly explicit information coming from the perception system of the ego car is entered in the knowledge graph, such as object measurements, infrastructural information, the classification of entities, and their relations. The perception module itself is out of the scope of this research. In the simulations, perfect sensors are used for the localization and classification of entities and lane markings. Although the risk-estimation method does incorporate stochastic properties of the measured quantities. For the classifier, the certainty of classification increases when the object is closer or is in the scene longer.  
Besides the explicit context information, we add implicit information from the interpretation of the scene in two different ways.  First, we calculate a collision risk between the ego and the perceived entities using the PRISMA~\cite{deGelder2023} method. And secondly, we use default reasoning rules to infer the risk levels of entities and the acceptability to cross lane markings. 
The classification, risk levels, and the acceptability to cross lane markings are queried from the knowledge graph and used to fill the risk field of the trajectory generator. In the next subsections, we elaborate on the method for calculating collision risk, the structure of the knowledge graph, and the situational rules.

\subsection{Collision Probability}
To calculate the collision probability of the ego vehicle with other objects in its lane, we employ the PRISMA method \cite{deGelder2023}. The probabilities are calculated analytically as described below. We start by calculating whether or not the vehicle is on a collision course with the nearest observed object edges. This is done based on the closest longitudinal and lateral distance to the $j$-th obstacle edge ($\bar{x}^{(j)}$,$\bar{y}^{(j)}$), the heading difference of the $j$-th obstacle edge ($\bar{\psi}^{(j)}$), the length of the  $j$-th obstacle edge ($w_t^{(j)}$), the width of the ego vehicle ($w_e$) and the longitudinal velocity of the ego vehicle ($v_k$). Note, that these measurements may involve a noise-induced uncertainty. Assuming that the ego vehicle does not change its course, we can determine the chance of collision. Let $\bar{y}^{(j)}_{\bar{x}^{(j)}=0}$ be the lateral distance to the obstacle if the ego vehicle continues to drive until $\bar{x}^{(j)}=0$. We can describe $\bar{y}^{(j)}_{\bar{x}^{(j)}=0}$ as:
\begin{align}
\bar{y}^{(j)}_{\bar{x}^{(j)}=0} = \bar{y}^{(j)} + \bar{x}^{(j)} \tan(\bar{\psi}^{(j)}).
\end{align}

This means that there is no collision when
\begin{align}
|\bar{y}^{(j)}_{\bar{x}^{(j)}=0}| > \frac{w_e + w_t^{(j)}}{2}\quad\forall j.
\end{align}

Assuming that $\bar{\psi}^{(j)}$ is small (such that $\tan(\bar{\psi}^{(j)})=\bar{\psi}^{(j)}$), we know the following:
\begin{align}
&V(\bar{y}^{(j)}_{\bar{x}^{(j)}=0}) = V(\bar{y}^{(j)}) + V(\bar{\psi}^{(j)})((\bar{x}^{(j)})^2 + V(\bar{x}^{(j)}))\\
&V\left( \frac{w_e + w_t^{(j)}}{2}\right) = \frac{1}{4}V\left(w_t^{(j)}\right)
\end{align}

Where the operation $V(\cdot)$ denotes the second moment, i.e., the variance. This means that the distance between the vehicle and the object at $\bar{x}^{(j)}=0$ is $|\bar{y}^{(j)}_{\bar{x}^{(j)}=0}|- \frac{w_t^{(j)}+w_e}{2}$ and has a variance of $V(\bar{y}^{(j)}_{\bar{x}^{(j)}=0})+\frac{1}{4}V(w_t^{(j)})$. If we assume that our variables are normally distributed ($\mathcal{N}(\mu,\,\bar{\sigma}^{2})$), we get the following probability of collision $P(c)$.

\begin{align}
P(c) &= P\left(-\frac{w_e + w_t^{(j)}}{2}<\mathcal{N}(\mu,\,\bar{\sigma}^{2})< \frac{w_e + w_t^{(j)}}{2}\right)\nonumber\\
\mu &= \bar{y}^{(j)}_{\bar{x}^{(j)}=0}\nonumber\\
\bar{\sigma}^{2} &= V(\bar{y}^{(j)}) + V(\bar{\psi}^{(j)})((\bar{x}^{(j)})^2 + V(\bar{x}^{(j)}))\nonumber
\end{align}

In addition, we will also add a term to the collision probability that considers the distance to the obstacle and the opportunity for the ego car to brake before arriving at $\bar{x}^{(j)}=0$. We can find the needed deceleration $d$ to arrive at $\bar{x}^{(j)}=0$ with $v_k=0$ with the formula:

\begin{align}
d = \frac{\left(v_k\cos({\bar{\psi}^{(j)}})\right)^2}{2\bar{x}^{(j)}}
\end{align}

It is then possible to set a threshold $d_0$ such that, if $d<d_0$, the ego vehicle has enough time to brake, and the collision probability is 0. But instead of using a hard threshold to decide if a collision is inevitable, the condition is smoothed by using the sigmoid function:

\begin{align}
S(u)\!=\!\frac{1}{1-e^{-u}}
,\:u = \alpha\left(d_0- \left(\!v_k\cos({\bar{\psi}^{(j)}})\right)^2{2\bar{x}^{(j)}}\!\right)
\end{align}

Where $\alpha$ is a parameter that determines the 'smoothness' of the condition. This makes the new probability calculation of a collision, $\bar{P}(c)$ with an object as follows:

\begin{align}
\bar{P}(c)=P(c) S\left(\alpha\left(d_0-\left(v_k\cos({\bar{\psi}^{(j)}})\right)^2{2\bar{x}^{(j)}}\right)\right)
\end{align}

\subsection{Knowledge Graph and Situational Rules}
\begin{figure}[t]
    \centering
\includegraphics[width=\columnwidth]{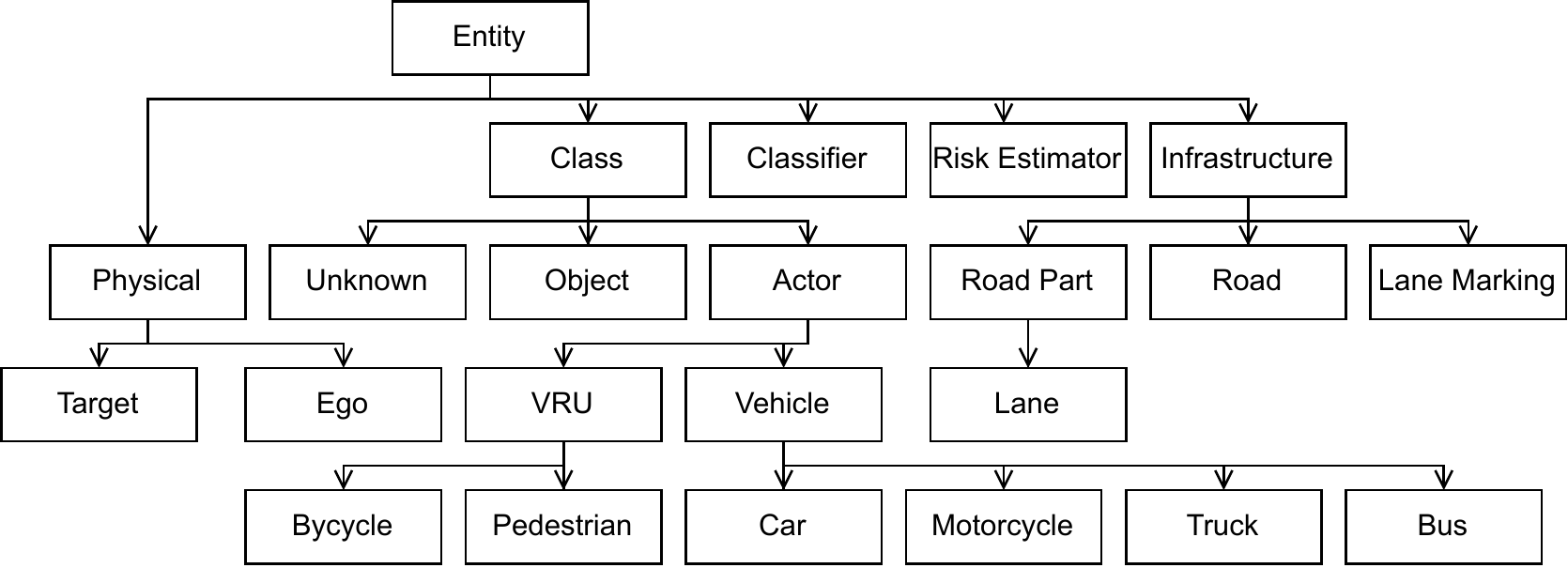}    
\caption{Entity hierarchy}
    \label{fig:entities}
    \vspace{-0.5cm}
\end{figure}

\begin{figure}[t]
    \centering
\includegraphics[width=0.6\columnwidth]{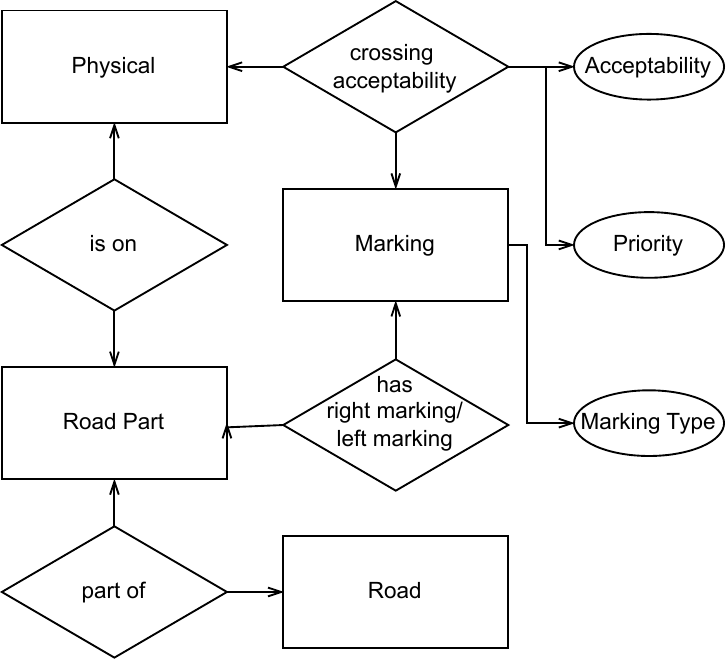}
    \caption{Most important relations for infrastructure entities.}
    \label{fig:infra_relations}
\end{figure}

\begin{figure}[t]
    \centering
\includegraphics[width=0.8\columnwidth]{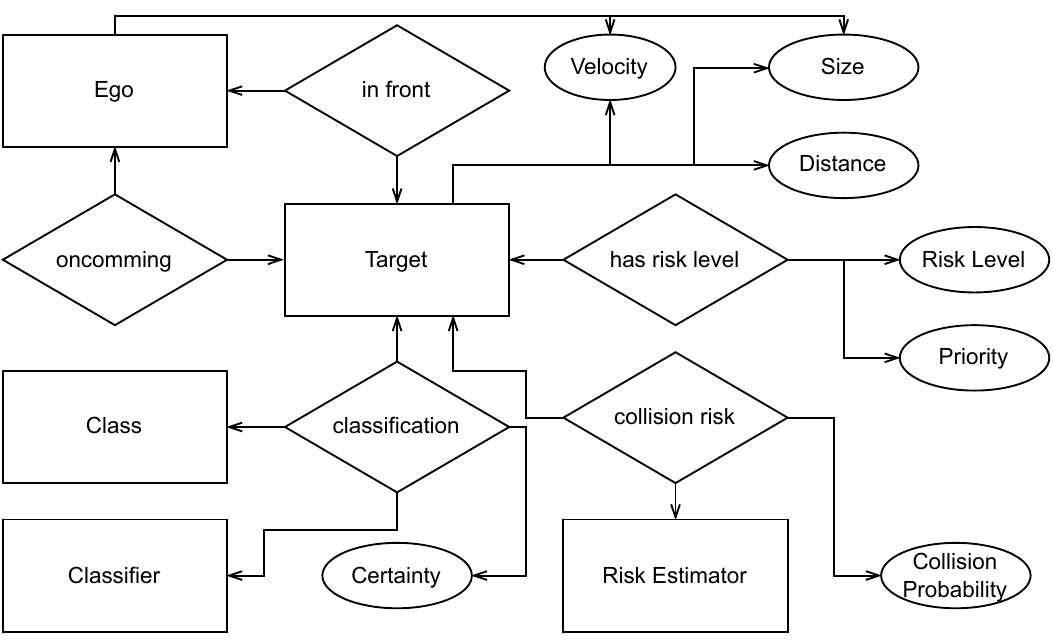}
    \caption{Most important relations for physical entities.}
    \label{fig:object_relations}
    \vspace{-0.5cm}
\end{figure}

The knowledge graph and the contextual rules are implemented using TypeDB \cite{pribadi_grakn_2020}. In the knowledge graph, the context is represented by a graph. There are three types of nodes in the graph: entities, relations, and attributes, from which all types of nodes are derived. There are two types of edges, which are the ‘role’ and ‘owns’ edge types. 
Entities are nodes that can exist in the graph independently from any other concept. An example of an entity can be a perceived object or a lane, but it can also be an object detector that is part of the automated driving system. Figure \ref{fig:entities} shows all the defined entities that are defined in the ontology of the used knowledge graph. Secondly, there are relations, which are nodes that are connected to other entities or relation nodes with role edges. An example of a relation is the ‘is\_on’ relation, which connects a physical entity and a road part entity together. Where the physical entity and road have their specific role connections to the ‘is\_on’ relation node. Attributes are nodes that have a specific data value connected to them such as a string, float, boolean, or integer.  Examples of attributes are longitudinal velocity, a classification certainty of an object, or an associated risk value. Attributes can be linked to entities and relations by an ‘owns’ or ‘role’ edge. 
Figure \ref{fig:infra_relations} shows the relations, entities, and attributes which are related to the infrastructure, and Figure \ref{fig:object_relations} shows the relations for the other road users that are detected.
One of the main advantages of using a knowledge graph is the reasoner, which is able to add extra implicit information to the situation using second-order predicate logic rules. More specifically, in TypeDB, rules are written as Horn clauses~\cite{horn_sentences_1951}, which can add relations to the graph. The implemented rules for our use cases add the ‘crossing acceptability’ to lane markings and the ‘has risk level’ to objects. The ‘crossing acceptability’ relation represents the implicit knowledge on whether it is acceptable to cross a lane marking for the ego based on rules for overtaking. The ‘has risk level’ relation represents the overall risk level of objects for the ego.
\begin{algorithm}[t]
     \caption{The rule which in the natural language would be 'If a lane marking is dashed it is acceptable to cross it'}
     \label{alg:example_rule}
\begin{algorithmic}
        \If{$\exists x,y | (x\in \text{lane}\_\text{marking})\land (y\in \text{ego}) \land (x.\text{owns}.\text{lane}\_\text{marking}\_\text{type} = \text{dashed})\qquad\qquad\qquad$}
            \State \textbf{insert} 
            \Indent
                \State $z | z \in \text{crossing}\_\text{acceptability}$
                \State $z.\text{lane}\_\text{marking} = x$
                \State $z.\text{ego} = y$
                \State $z.\text{acceptability} = 1$
                \State $z.\text{priority} = 0$
            \EndIndent
        \EndIf
\end{algorithmic}
\end{algorithm}
These kind of rules, however, rely on common sense reasoning and have to deal with exceptions, conflicting rules, and an open world. Therefore these rules are implemented as prioritized default rules and do not only link a risk or acceptability attribute to a lane marking or object, but also add a priority. An example of such a rule is shown in Alg.~\ref{alg:example_rule}. When extracting the actual acceptability and risk levels, only the relation with the highest priority is considered. 
The rules that add a ‘crossing acceptability’ relation and their priorities are listed in Table~\ref{table:line rules}. They are based on common sense and traffic rules on how to interact when wanting to overtake on an urban two-way road. The rules that add a ‘has risk level’ consider three possible risk levels: low, medium, and high. The rules and their priorities are listed in Table~\ref{table:risk rules}.
\begin{table}[t]
\begin{center}
\begin{tabular}{ |c | p{6.5cm} |}
 \hline
 \textbf{Priority} & \textbf{Natural language rule} \\
  \hline
 0 & When a line is dashed, it is acceptable to cross \\
  \hline
 0 & When a line is solid, it is not acceptable to cross\\
  \hline
 1 & When a line is dashed and it is on the left side of your lane, it is not acceptable to cross\\
  \hline
2& When a line is dashed and it is on the left side of your lane and there is a object in front of you that is closer than 20 meters, the line is acceptable to cross\\
 \hline
  3 & When a line is dashed and it is on the left side of your lane and there is an oncoming vehicle on the lane left of you closer than 50 meters, the line is not acceptable to cross\\
   \hline
4 & When a line is dashed and there is a VRU crossing the road, it is not acceptable to cross the line.\\
 \hline
\end{tabular}
\caption{Natural language versions of the rules about crossing acceptability.}
\label{table:line rules}
\end{center}
\vspace{-0.7cm}
\end{table}
\begin{table}[t]
\begin{center}
\begin{tabular}{ |c | p{6.5cm} |}
\hline
 \textbf{Priority} & \textbf{Natural language rule} \\
  \hline
 0 & A object has a medium risk. \\  \hline
 1 & If a object has a classification of an artificial object with a certainty above 0.8 and is smaller than 0.4 m on all dimensions, the object's risk is low. \\  \hline
 2 & If a object has a classification as a VRU with a certainty value above 0.05 and has a collision probability higher than 0.05, the object's risk is high\\  \hline
 3 & If a object has a collision probability of higher than 0.2, the object's risk is high\\
 \hline
\end{tabular}
\caption{Natural language versions of the rules about object risk levels.}
\label{table:risk rules}
\end{center}
\vspace{-0.5cm}
\end{table}
The conclusions of these prioritized rules are used by the trajectory generator. As stated, the rules themselves are not directly converted into decisions. Instead, the 'crossing acceptability' relation of lane markings and the 'risk level' relation of objects are used to decide how the risk potential field is populated. For each lane marking there is a highest prioritized crossing acceptability relation which is either 1 or 0. Based on the type of line and the acceptability an amplitude of the potential field is selected. This is done in Table~\ref{table:lines weights}. For the objects, the amplitude and standard deviation of their risk fields are dependent on the classification, size of the object, and risk level as provided in Table~\ref{table:objects weights}. 

\begin{table}[t]
\begin{center}
\begin{tabular}{ p{1.5cm}  p{1.5cm} p{1cm} p{1cm}}
 \toprule
 \textbf{Line type} & \textbf{acceptability} & $A_L$ & $\sigma$\\
 \toprule
solid & 0 & 4 & 0.6\\
solid & 1 & 0 & 0.6\\
dashed & 0 & 1.5 & 0.6\\ 
dashed & 1 & 0 & 0.6\\
 \toprule
\end{tabular}
\caption{Risk parameters for markings based on type and crossing acceptability}
\label{table:lines weights}
\end{center}
\vspace{-0.5cm}
\end{table}
\begin{table}[t]
\begin{center}
\begin{tabular}{ p{1.5cm}  p{1.1cm} p{0.5cm} p{1.5cm} p{1.6cm} }
 \toprule
 \textbf{Object class} & \textbf{Risk} & $A_O$ & $\sigma_x$ & $\sigma_y$\\
 \toprule
Car & low & 2 & $1.5l_O+0.05$ & $1.2w_O+0.05$ \\
Car & medium & 3 & $1.5l_O+0.1$ & $1.2w_O+0.1$ \\
Car & high & 4 & $1.5l_O+0.3$ &$ 1.2w_O+0.3$ \\
Art. object & low & 1 & $1.5l_O+0.2$ &$ 1.2w_O+0.2$ \\
Art. object & medium & 2 & $1.5l_O+0.25$ & $1.2w_O+0.25$ \\
Art. object & high & 3 & $1.5l_O+0.4$ & $1.2w_O+0.4$ \\
Pedestrian & low & 2 & $1.5l_O+1.2 $& $1.2w_O+1.2$ \\
Pedestrian & medium & 3 & $1.5l_O+1.7$ &$ 1.2w_O+1.7$ \\
Pedestrian & high & 4 &$ 1.5l_O+2.2$ & $1.2w_O+2.2$ \\
 \toprule
\end{tabular}
\caption{Risk parameters for objects}
\label{table:objects weights}
\end{center}
\vspace{-0.5cm}
\end{table}

\section{Situation-aware motion planning}~\label{sec:method1}
When populating the problem space of the MPC trajectory generator, information from the knowledge graph is collected. Based on the classification of objects and the inferred risk, the amplitude and standard deviation of risk potential fields around objects are determined. The risk field represents an adaptive version of the risk field proposed in~\cite[Eq. 6]{van_der_ploeg_long_2022}, and is represented as follows.
\begin{align}
    \widetilde{\mathcal{U}}^{O}_k=&\sum_{n=1}^{N_o}A_O^{(n)}e^{-\frac{f_o^{(n)}(\mathbf{x}_k,\mathbf{o}_k)}{2}},\forall n\in \mathbb{Z}^+\label{eq:objects}\\
    f_o^{(n)}=&\begin{bmatrix}x_k-\widehat{x}^{(n)}_k& y_k-\widehat{y}^{(n)}_k\end{bmatrix}R\Sigma^{-1} R^\intercal\begin{bmatrix}x_k-\widehat{x}^{(n)}_k\\ y_k-\widehat{y}^{(n)}_k\end{bmatrix}\nonumber\\
    \Sigma=&\begin{bmatrix}{(\sigma_x^{(n)})^2}&0\\0&{(\sigma_y^{(n)})^2}\end{bmatrix},\quad R=\begin{bmatrix}\cos\widehat{\theta}^{(n)}_k&-\sin\widehat{\theta}^{(n)}_k\\\sin\widehat{\theta}^{(n)}_k&\cos\widehat{\theta}^{(n)}_k\end{bmatrix}\nonumber
\end{align}
$(\widehat{x}_k^{(n)},\widehat{y}_k^{(n)},\widehat{\theta}_k^{(n)})$ represents the state of the $n$-th object, $N_o$ represents the number of objects. The adaptive terms $A_O^{(n)},\sigma_x^{(n)},\sigma_y^{(n)}$ represent the amplitude and standard deviations of the bi-variate risk-field of the $n$-th object, as derived in Sec.~\ref{sec:method0}.
Similarly, for the road markings, the adaptive risk field proposed is represented as follows.
 \begin{align}
    \label{eq:infra}
    \widetilde{\mathcal{U}}^{I}_k=&\sum_{i=1}^{N_d}A_I^{(i)}e^{-\frac{\left(y_k-\tilde{y}_k^{(i)}(\tilde{x}^{(i)}_k)\right)^2}{2 (\sigma^{(i)})^2}},\forall i\in \mathbb{Z}^+
\end{align}
where $N_d$ represents the number of road lines and the adaptive terms $A_I^{(i)}$, $\sigma^{(i)}$ represents the amplitude of the risk field and the standard deviation of the risk field of the infrastructure, respectively, as derived in Sec.~\ref{sec:method0}. The adaptive formulation of the risk-fields can then be implemented into the model-predictive planning problem as introduced in our previous work~\cite{van_der_ploeg_long_2022}, as follows
 \begin{subequations}
\label{eq:mpcprogram}
    \begin{align}
        J^\star_{0\rightarrow N}(\mathbf{x}_s)=&\min_{u,x}\sum_{k=1}^{N-1}\ell_\mathbf{x}\left(\mathbf{x}_k-\mathbf{r}_k\right)+\sum_{k=0}^{N-1}\ell_\mathbf{u}\left(\mathbf{u}_k\right)\nonumber\\+m&(\mathbf{x}_{N}-\mathbf{r}_{N})+\widetilde{\mathcal{U}}_k^O+\widetilde{\mathcal{U}}_k^I\label{eq:pathmpc1}\\\text{s.t.}\quad\mathbf{x}_0=&\mathbf{x}_s\label{eq:pathmpc2}\\
        \quad \mathbf{x}_{k+1}=&f(\mathbf{x}_k,\mathbf{u}_k),\quad\forall k\geq 0\label{eq:pathmpc3}\\
        \quad g(\mathbf{x}_k,\mathbf{u}_k)\leq&\: 0,\quad \forall k\geq 0\label{eq:pathmpc4}
    \end{align}
\end{subequations}
where, in~\eqref{eq:pathmpc1}, the functions $\ell_{\mathbf{x}}(\cdot),\:\ell_{\mathbf{u}}(\cdot)$ quadratically penalize the stage error $\mathbf{x}_k-\mathbf{r}_k$, and the control input $\mathbf{u}_k$, respectively. Furthermore, the term $m(\cdot)$ quadratically penalizes the terminal state at time horizon $N$. The state $\mathbf{x}_s$ in~\eqref{eq:pathmpc2} represents the initial state. The function $f$ in~\eqref{eq:pathmpc3} represents the non-linear dynamics as depicted in~\eqref{eq:vehiclemodel}. Finally, the constraints in~\eqref{eq:pathmpc4} enforce the boundedness of the state $\mathbf{x}$ and the input $\mathbf{u}$ through enforcement of the sets~\eqref{eq:boundedsets}. 
\section{Simulation study}~\label{sec:simulation}
\begin{figure}[t]
    \centering
     \fontsize{9pt}{11pt}\selectfont
    \def\svgwidth{1\columnwidth}
\begingroup%
  \makeatletter%
  \providecommand\color[2][]{%
    \errmessage{(Inkscape) Color is used for the text in Inkscape, but the package 'color.sty' is not loaded}%
    \renewcommand\color[2][]{}%
  }%
  \providecommand\transparent[1]{%
    \errmessage{(Inkscape) Transparency is used (non-zero) for the text in Inkscape, but the package 'transparent.sty' is not loaded}%
    \renewcommand\transparent[1]{}%
  }%
  \providecommand\rotatebox[2]{#2}%
  \newcommand*\fsize{\dimexpr\f@size pt\relax}%
  \newcommand*\lineheight[1]{\fontsize{\fsize}{#1\fsize}\selectfont}%
  \ifx\svgwidth\undefined%
    \setlength{\unitlength}{1995bp}%
    \ifx\svgscale\undefined%
      \relax%
    \else%
      \setlength{\unitlength}{\unitlength * \real{\svgscale}}%
    \fi%
  \else%
    \setlength{\unitlength}{\svgwidth}%
  \fi%
  \global\let\svgwidth\undefined%
  \global\let\svgscale\undefined%
  \makeatother%
  \begin{picture}(1,0.73796992)%
    \lineheight{1}%
    \setlength\tabcolsep{0pt}%
    \put(0,0){\includegraphics[width=\unitlength,page=1]{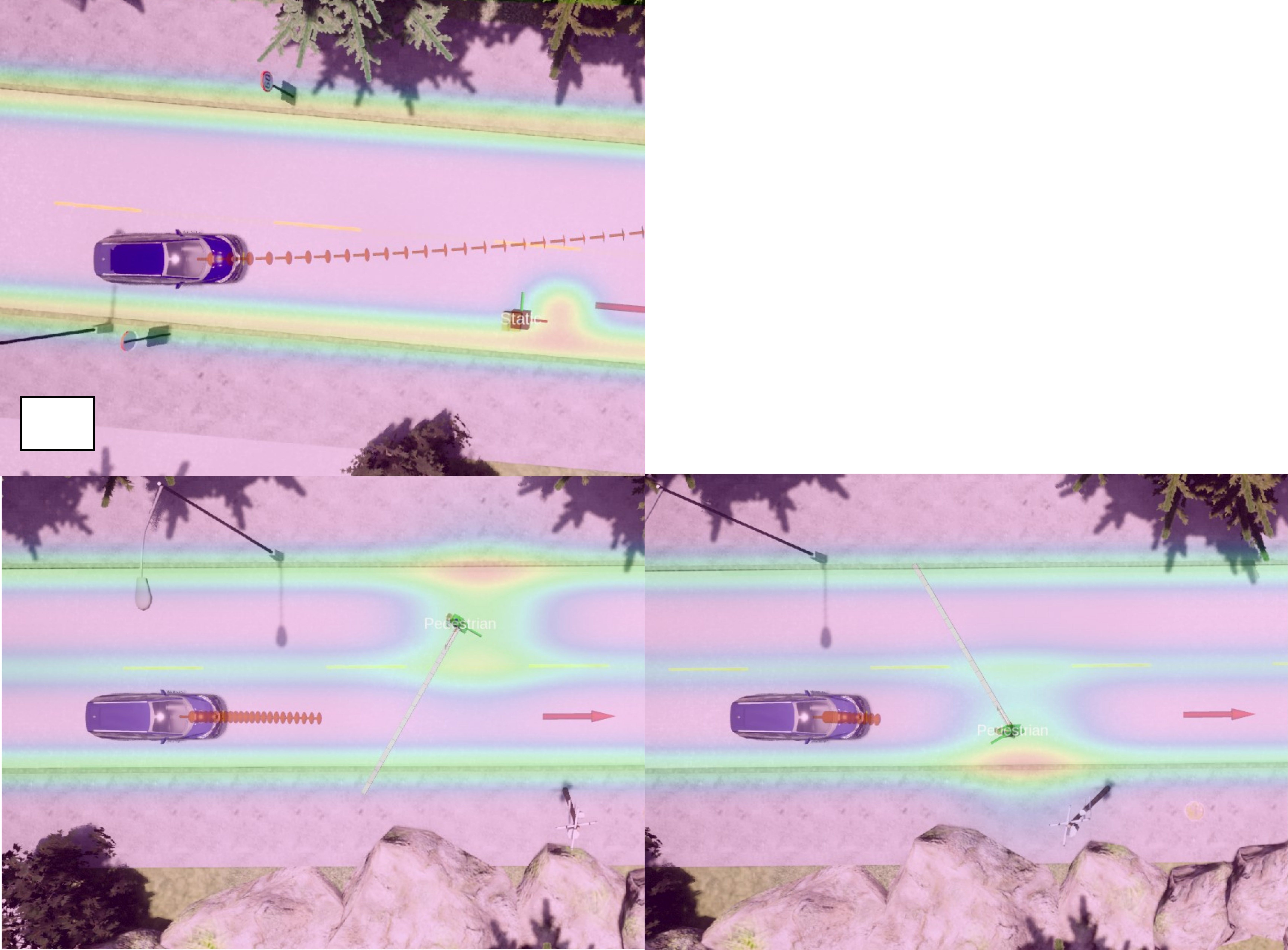}}%
    \put(0.0291447,0.39586466){\makebox(0,0)[lt]{\lineheight{1.25}\smash{\begin{tabular}[t]{l}a.\end{tabular}}}}%
    \put(0,0){\includegraphics[width=\unitlength,page=2]{Carla_scenarios.pdf}}%
    \put(0.0291447,0.02593985){\makebox(0,0)[lt]{\lineheight{1.25}\smash{\begin{tabular}[t]{l}c.\end{tabular}}}}%
    \put(0,0){\includegraphics[width=\unitlength,page=3]{Carla_scenarios.pdf}}%
    \put(0.52934771,0.02593985){\makebox(0,0)[lt]{\lineheight{1.25}\smash{\begin{tabular}[t]{l}d.\end{tabular}}}}%
    \put(0,0){\includegraphics[width=\unitlength,page=4]{Carla_scenarios.pdf}}%
    \put(0.52934771,0.39586466){\makebox(0,0)[lt]{\lineheight{1.25}\smash{\begin{tabular}[t]{l}b.\end{tabular}}}}%
  \end{picture}%
\endgroup%

    \caption{Depiction of the different use-cases in CARLA.}
    \label{fig:carlascenarios}
    \vspace{-0.5cm}
\end{figure}

In order to demonstrate the effectiveness of our proposed approach, the model-predictive planner, incorporated with the knowledge graph, is implemented in ROS and simulated in CARLA~\cite{dosovitskiy_carla_2017}. The model-predictive planner is programmed using the CasADi~\cite{andersson_casadi_2019} optimization framework, in combination with the IPOPT~\cite{wachter_implementation_2006} non-linear solver. For enabling these simulations, the object observations are directly provided to the vehicle, which makes a constant velocity-based prediction of the observations over the prediction horizon of the model-predictive planner. For longitudinal and lateral vehicle actuation, two feedback controllers are used which are proportional and proportional-derivative, respectively, identical to the controller used in~\cite[Sec. IV]{smit_informed_2022}. The vehicle is always initiated from the same position from standstill. In this section, four scenarios are demonstrated as depicted in Fig.~\ref{fig:carlascenarios} and described as follows.
\begin{enumerate}[label=\alph*.]
\item The ego vehicle is driving towards a static obstacle . (no. of variations: $9$)
\item The ego vehicle is driving towards a static obstacle with opposing traffic. (no. of variations: $108$)
\item The ego vehicle is driving towards a pedestrian who is crossing the road from the left. (no. of variations: $96$)
\item The ego vehicle is driving towards a pedestrian who is crossing the road from the right. (no. of variations: $96$)
\end{enumerate}
The actors and elements in these scenarios are varied in multiple degrees of freedom to show the robustness of our proposed approach. In scenario a., the initial position of the object in the ego lane is varied and the type of object is chosen from three categories: a vehicle, a glass bin, and a cardboard box. In scenario b., the initial position and velocity of the oncoming vehicle are varied, where the velocity of the oncoming vehicle is selected between $5.5-11\:\text{m}\cdot\text{s}^{-1}$. The variations of scenario b. are identical to a., although in the presence of oncoming traffic. In scenario c. and d., the initial position and heading of the pedestrian is varied, and its velocity ranges from $1-2.2\:\text{m}\cdot\text{s}^{-1}$.
\begin{table}
\small
\begin{tabular}{p{0.3cm} p{3.4cm} p{2.5cm} p{0.6cm}}\toprule
\textbf{Par.}              & \textbf{Description} & \textbf{Value} & \textbf{Unit}\\
\toprule
$t_s$ & Sampling time & 0.15 & s\\
$N$            & Trajectory horizon & 23 & -\\
$\ell_\mathbf{u}$           & Input weight & $[0.01,1]$ & -\\
$\ell_\mathbf{x}$ & Stage weight & $[0,0,0,0.001,0]$& -\\
$m$&Terminal weight&$[20,2,0,1,0]\cdot 10^{-2}$&-\\
$\underline{a},\overline{a}$&Long. acceleration bounds&$[-4,2]$&$\text{ms}^{-2}$\\
$\underline{\delta},\overline{\delta}$&Steering angle bounds&$[-0.1,0.1]$&$\text{rad}$\\
$\underline{v},\overline{v}$&Long. velocity bounds&$[-2,10]$&$\text{ms}^{-1}$\\
$\underline{\omega},\overline{w}$&Steering rate bounds &$[-1.5,1.5]$&$\text{rads}^{-1}$\\
\toprule
\end{tabular}
\caption{Simulation parameters}
\label{table:pars}
\vspace{-0.5cm}
\end{table}
The planner and knowledge-graph parameters (Tab.~\ref{table:objects weights},~\ref{table:lines weights},~\ref{table:pars}) remain constant throughout all simulations, showing the strength of the novel combination. For each scenario, we analyze the performance indicators of interest, e.g., no. of collisions, reaching the area of interest, minimal safety distances, and various comfort parameters.
\subsection{Object without oncoming traffic}
\begin{figure}
    \centering
    \includegraphics[width=\columnwidth]{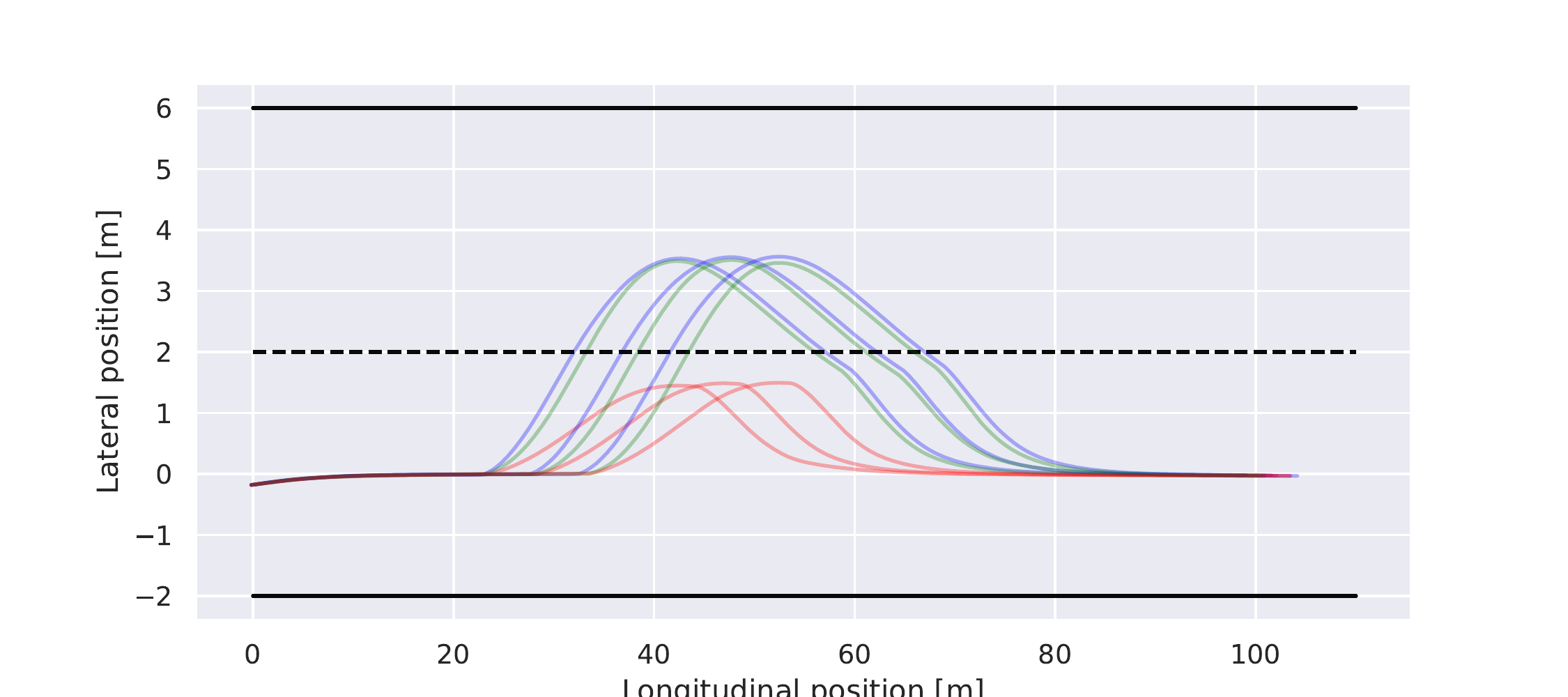}
    \caption{Driven trajectories for the object scenario without oncoming traffic.}
    \label{fig:res_traj_obj_without}
    \vspace{-0.5cm}
\end{figure}

\begin{figure}
    \centering
    \includegraphics[width=\columnwidth]{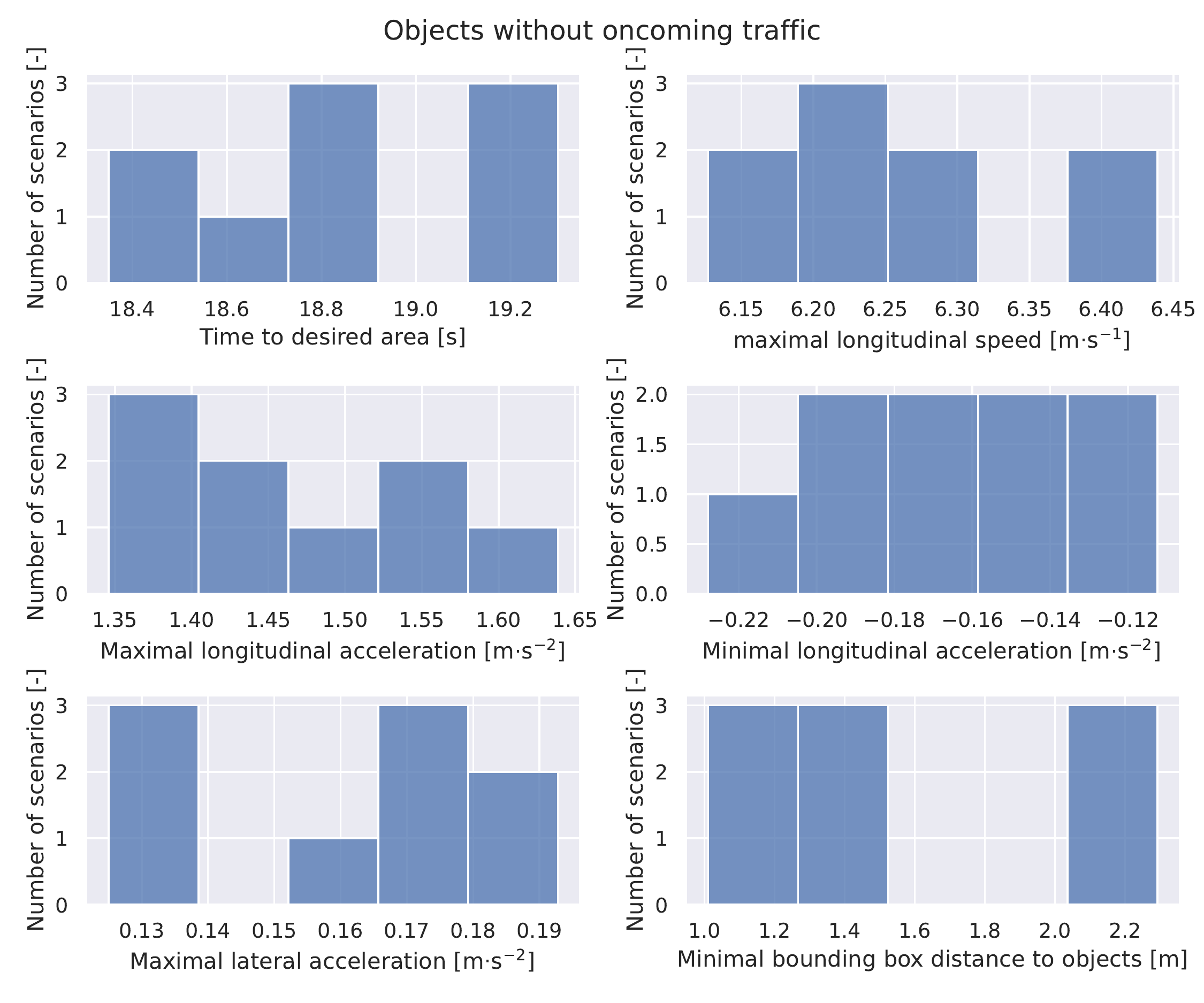}
    \caption{Results for the object scenario w/o oncoming traffic.}
    \label{fig:res_obj_without}
    \vspace{-0.6cm}
\end{figure}

 Out of the 9 simulated variations, three different object types with three different longitudinal initial positions of the objects are used. Since the vehicle starts from a standstill, the variations in initial position will result in different speeds with which the vehicle approaches the object. No crashes were recorded and the vehicle had always been able to reach its destination. The driven trajectories for these scenarios can be observed in Fig.~\ref{fig:res_traj_obj_without}. Here, the trajectories are color-labeled, where red depicts the trajectories with a cardboard box on the road, blue depicts the trajectories with a vehicle parked on the side of the road and green depicts the trajectories with a glass bin on the side of the road. In Fig.~\ref{fig:res_obj_without}, it can be seen that the time it takes to reach the destination is all within a second from each other, between $18.4-19.4$s. Moreover, the longitudinal speed and maximum/minimum acceleration is well within the MPC bounds. The lateral acceleration is within $0.2\text{m}\cdot \text{s}^{-2}$, which resides in the most comfortable region of lateral acceleration perceived by drivers~\cite{de_winkel_standards_2023}. Finally, the bounding box distance to the objects shows the consistency of the approach, where each object type is allocated a bin (the box is evaded with a margin of $1$m, the glass bin with $1.5$m, the vehicle with $2.25$m).
\subsection{Object with oncoming traffic}
\begin{figure}
    \centering
    \includegraphics[width=\columnwidth]{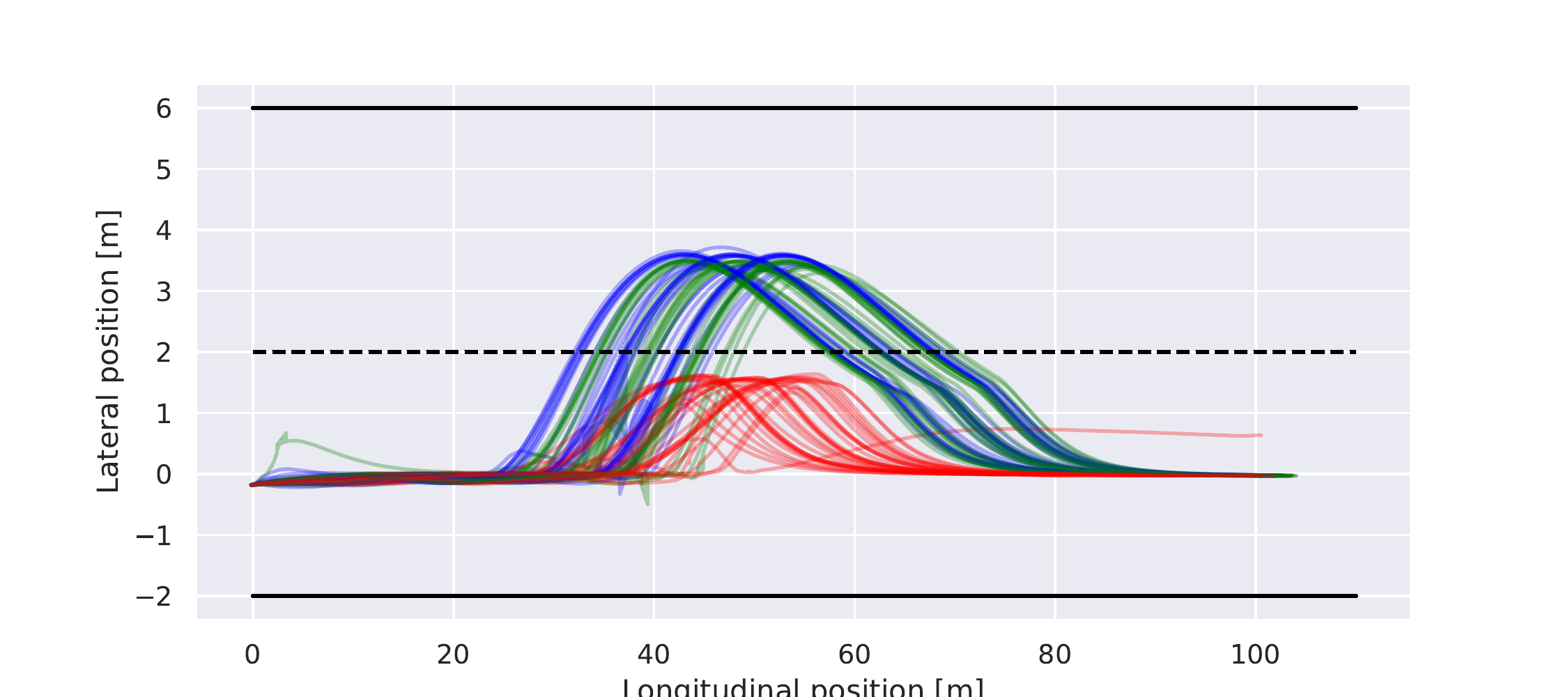}
    \caption{Driven trajectories for the object scenario with oncoming traffic.}
    \label{fig:res_traj_obj}
    \vspace{-0.5cm}
\end{figure}

\begin{figure}
    \centering
    \includegraphics[width=\columnwidth]{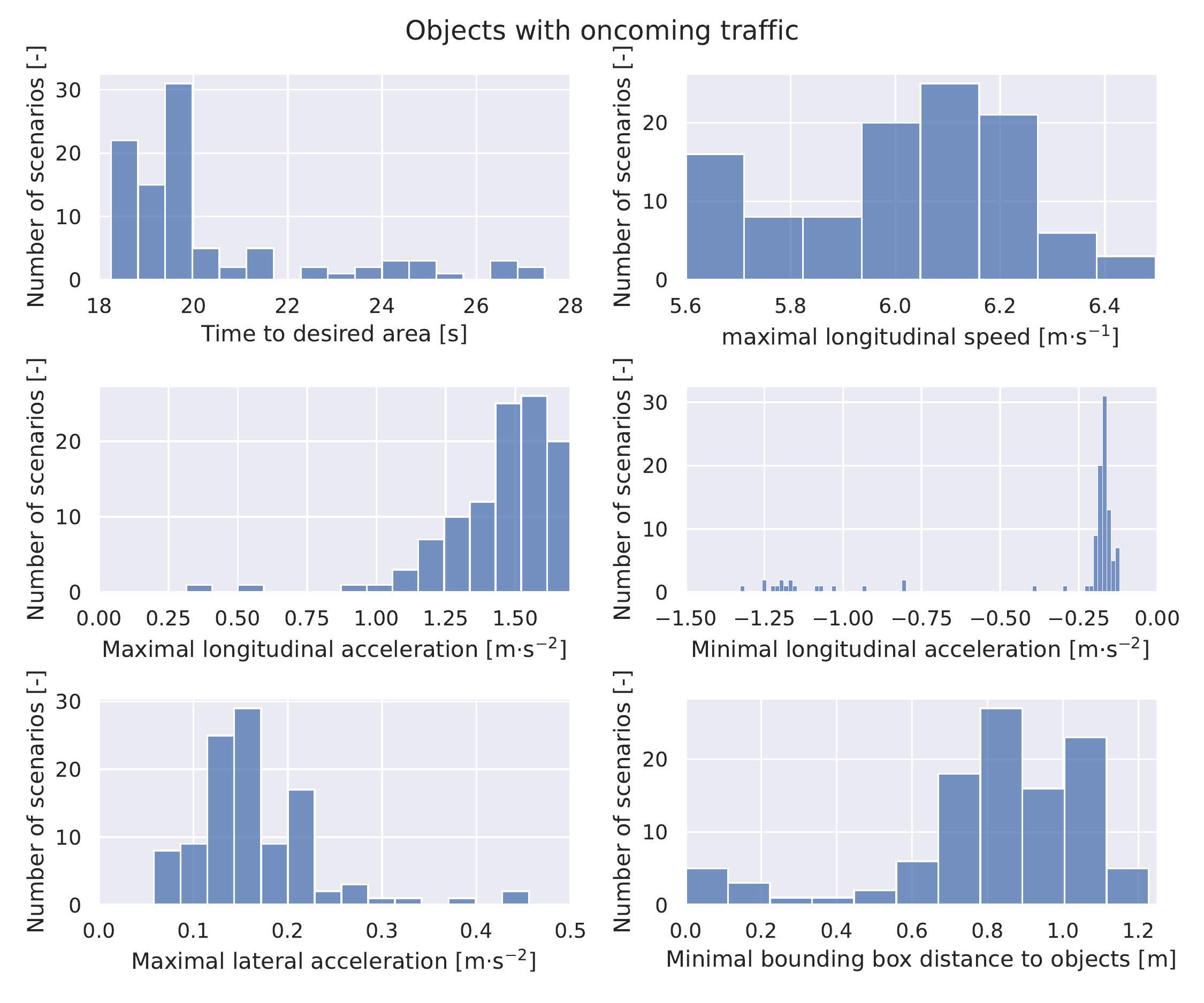}
    \caption{Results for the object scenario with oncoming traffic.}
    \label{fig:res_obj}
    \vspace{-0.5cm}
\end{figure}

 Out of the 108 simulated variations, one time the vehicle collided against an object (a cardboard box), and two times it failed to reach its destination. When replaying the scenarios, the root cause of the vehicle hitting the box, or not reaching its destination can be explained by the fact that it, during a lateral maneuver, had to turn back due to opposing traffic. During the collision scenario with the cardboard box (also seen as the red outlier in Fig.~\ref{fig:res_traj_obj}), the vehicle had to do a last-minute lane-change due to late detected oncoming traffic, as a result, the vehicle collided against the box and dragged it along, resulting in off-center driving. Although this is a low-risk maneuver, ideally we would strive for zero collisions. This case could potentially be solved by tuning the steering rate and deceleration limits, although this may come at the cost of comfort. When the vehicle is too close to the object, it could end up in a "dead-lock", since the cost associated with driving rearward and overtaking the object could be in perfect equilibrium with the cost associated with reaching the goal point, as a result, the vehicle would stand still during the remainder of the simulation. Although rare, this case could be circumvented by placing less or no weight on tracking the desired velocity when within a certain proximity to an object. The visual grouping of scenarios is again apparent. When looking at all other signals of interest in Fig.~\ref{fig:res_obj}, most of the scenarios are finished within a time span of $20s$, all other scenarios and, particularly the outliers above $26s$ can be explained by the fact that the vehicle has to wait for the oncoming traffic to pass and, potentially, drive rearward to be able to steer away from the object in front. The lateral acceleration is similar to the case without oncoming vehicles, although with a few outliers caused by late detections of the oncoming traffic, hence forcing the vehicle to quickly transition back to the original lane. The longitudinal velocity and maximum/minimum longitudinal acceleration are well within the bounds of the program, the lateral acceleration is mostly within $0.3\text{m}\cdot\text{s}^{-2}$, which resides in the most comfortable region of lateral acceleration perceived by drivers~\cite{de_winkel_standards_2023}. Finally, sufficient distance is kept between the bounding boxes, as the smallest values between $0-0.2$m are associated with overtaking the cardboard box, which is a very low-risk maneuver. 
\subsection{Pedestrian crossing from left and right}
\begin{figure}
    \centering
    \includegraphics[width=\columnwidth]{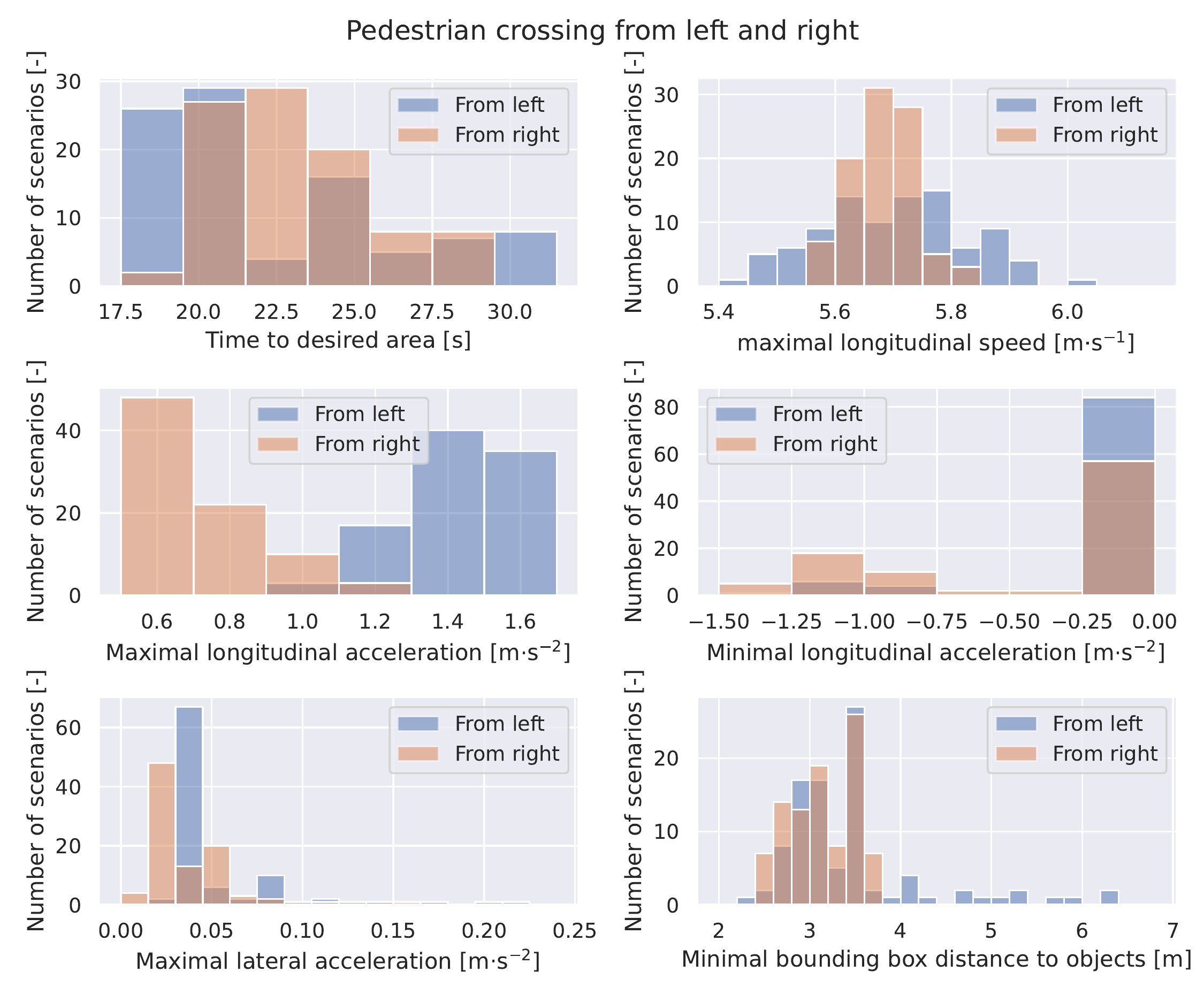}
    \caption{Results for the pedestrian crossing scenario.}
    \label{fig:res_ped_all}
    \vspace{-0.5cm}
\end{figure}

The simulated data for the pedestrian crossing scenarios from the left and right (each consisting of $96$ variations) are provided in Fig.~\ref{fig:res_ped_all}. In all variations, the vehicle arrives at the desired area without collision. The trajectories are omitted from the paper, as all trajectories follow the lane, and, as such, little lateral action can visually be observed. This can be explained due to the incorporated social norm that it is not appropriate to overtake a pedestrian when it is crossing the road. As such, the vehicle will brake with a sufficient standstill distance (above $2.5$m for all simulated variants) and wait for the pedestrian to pass, after which it will proceed to drive to its destination. The vehicle anticipates well, as it brakes with a low deceleration (mostly within $-0.25\text{m}\cdot\text{s}^{-2}$), and little lateral acceleration is observed since the vehicle stays within its lane. When accelerating, the maximum acceleration is well within the MPC bounds.
\section{Conclusion}~\label{sec:conclusion}
In this work, we have proposed a novel combination of a model-predictive generator with a knowledge graph to incorporate contextual information. As opposed to common approaches in the literature, no discrete decisions about the type of maneuvers are made. Instead, the trajectory generator only uses the risk field, which is composed of explicit observations, traffic rules, and social norms, as well as dynamic constraints to generate trajectories to an area of interest. Our method is demonstrated through $309$ simulation variations, where it has been shown to be robust and effective. Future work includes: preventing "dead-lock" situations as were seen in the simulation results, and the incorporation of prediction models for seen and unseen/occluded objects.
    \printbibliography

@INPROCEEDINGS{Huang_2019,
  author={Huang, Lu and Liang, Huawei and Yu, Biao and Li, Bichun and Zhu, Hui},
  booktitle={2019 4th Asia-Pacific Conference on Intelligent Robot Systems (ACIRS)}, 
  title={Ontology-Based Driving Scene Modeling, Situation Assessment and Decision Making for Autonomous Vehicles}, 
  year={2019},
  volume={},
  number={},
  pages={57-62},
  doi={10.1109/ACIRS.2019.8935984}}

@article{Zhao_2017,
author = {Zhao, Lihua and Ichise, Ryutaro and Liu, Zheng and Mita, Seiichi and Sasaki, Yutaka},
year = {2017},
month = {07},
pages = {1425-1439},
title = {Ontology-Based Driving Decision Making: A Feasibility Study at Uncontrolled Intersections},
volume = {E100.D},
journal = {IEICE Transactions on Information and Systems},
doi = {10.1587/transinf.2016EDP7337}
}

@preamble{ " \newcommand{\noop}[1]{} " }

@article{deGelder2023,
  author = {de Gelder, E. and Adjenughwure, K. and Manders, J. and Snijders, R. and Paardekooper, J.-P. and Op den Camp, O.  and Tejada Ruiz, A. and De Schutter, B.}, 
  title   = {{PRISMA}: A Novel Approach for Deriving Probabilistic Surrogate Safety Measures for Risk Evaluation},
  year    = {\noop{2023} Under Review},
  }

@article{wachter_implementation_2006,
	title = {On the implementation of an interior-point filter line-search algorithm for large-scale nonlinear programming},
	volume = {106},
	issn = {1436-4646},
	doi = {10.1007/s10107-004-0559-y},
	number = {1},
	journal = {Mathematical Programming},
	author = {Wächter, A. and Biegler, L.T.},
	year = {2006},
	pages = {25--57},
}

@inproceedings{halilaj_knowledge_2021,
	title = {A {Knowledge} {Graph}-{Based} {Approach} for {Situation} {Comprehension} in {Driving} {Scenarios}},
	booktitle = {Lecture {Notes} in {Computer} {Science}},
	author = {Halilaj, Lavdim and Dindorkar, Ishan and Lüttin, Jürgen and Rothermel, Susanne},
	year = {2021},
	pages = {699--716},
}

@article{andersson_casadi_2019,
	title = {{CasADi}: a software framework for nonlinear optimization and optimal control},
	volume = {11},
	issn = {1867-2949},
	shorttitle = {{CasADi}},
	doi = {10.1007/s12532-018-0139-4},
	journal = {Mathematical Programming Computation},
	author = {Andersson, J.A.E. and Gillis, J. and Horn, G. and Rawlings, J.B. and Diehl, M.},
	year = {2019},
	pages = {1--36},
}

@article{horn_sentences_1951,
	title = {On {Sentences} {Which} are {True} of {Direct} {Unions} of {Algebras}},
	volume = {16},
	issn = {0022-4812},
	url = {https://www.jstor.org/stable/2268661},
	doi = {10.2307/2268661},
	number = {1},
	urldate = {2022-11-25},
	journal = {The Journal of Symbolic Logic},
	author = {Horn, Alfred},
	year = {1951},
	pages = {14--21},
}

@article{gutjahr_lateral_2017,
	title = {Lateral {Vehicle} {Trajectory} {Optimization} {Using} {Constrained} {Linear} {Time}-{Varying} {MPC}},
	volume = {18},
	issn = {1558-0016},
	doi = {10.1109/TITS.2016.2614705},
	number = {6},
	journal = {IEEE Transactions on Intelligent Transportation Systems},
	author = {Gutjahr, Benjamin and Gröll, Lutz and Werling, Moritz},
	month = jun,
	year = {2017},
	pages = {1586--1595},
}

@article{xie_distributed_2022,
	title = {Distributed {Motion} {Planning} for {Safe} {Autonomous} {Vehicle} {Overtaking} via {Artificial} {Potential} {Field}},
	volume = {23},
	issn = {1558-0016},
	doi = {10.1109/TITS.2022.3189741},
	number = {11},
	journal = {IEEE Transactions on Intelligent Transportation Systems},
	author = {Xie, Songtao and Hu, Junyan and Bhowmick, Parijat and Ding, Zhengtao and Arvin, Farshad},
	month = nov,
	year = {2022},
	pages = {21531--21547},
}

@article{ji_path_2017,
	title = {Path {Planning} and {Tracking} for {Vehicle} {Collision} {Avoidance} {Based} on {Model} {Predictive} {Control} {With} {Multiconstraints}},
	volume = {66},
	issn = {1939-9359},
	doi = {10.1109/TVT.2016.2555853},
	number = {2},
	journal = {IEEE Transactions on Vehicular Technology},
	author = {Ji, Jie and Khajepour, Amir and Melek, Wael William and Huang, Yanjun},
	month = feb,
	year = {2017},
	pages = {952--964},
}

@inproceedings{chen_socially_2017,
	title = {Socially aware motion planning with deep reinforcement learning},
	volume = {2017-September},
	isbn = {978-1-5386-2682-5},
	doi = {10.1109/IROS.2017.8202312},
	author = {Chen, Y.F. and Everett, M. and Liu, M. and How, J.P.},
	year = {2017},
	pages = {1343--1350},
}

@article{paardekooper_hybrid-ai_2021,
	title = {A hybrid-{AI} approach for competence assessment of automated driving functions},
	url = {https://repository.ubn.ru.nl/handle/2066/246424},
	urldate = {2023-01-18},
	journal = {SafeAI@ AAAI},
	author = {Paardekooper, J. P. and Comi, M. and Grappiolo, C. and Snijders, R. R. and Vught, W. van and Beekelaar, R.},
	year = {2021},
}

@article{endsley_toward_1995,
	title = {Toward a {Theory} of {Situation} {Awareness} in {Dynamic} {Systems}},
	volume = {37},
	issn = {0018-7208},
	url = {https://doi.org/10.1518/001872095779049543},
	doi = {10.1518/001872095779049543},
	language = {en},
	number = {1},
	urldate = {2023-01-20},
	journal = {Human Factors},
	author = {Endsley, Mica R.},
	month = mar,
	year = {1995},
	pages = {32--64},
}

@article{hang_human-like_2021,
	title = {Human-{Like} {Decision} {Making} for {Autonomous} {Driving}: {A} {Noncooperative} {Game} {Theoretic} {Approach}},
	volume = {22},
	issn = {1524-9050},
	shorttitle = {Human-{Like} {Decision} {Making} for {Autonomous} {Driving}},
	doi = {10.1109/TITS.2020.3036984},
	language = {English},
	number = {4},
	journal = {IEEE Transactions on Intelligent Transportation Systems},
	author = {Hang, P. and Lv, C. and Xing, Y. and Huang, C. and Hu, Z.},
	year = {2021},
	pages = {2076--2087},
}

@article{de_winkel_standards_2023,
	title = {Standards for passenger comfort in automated vehicles: {Acceleration} and jerk},
	volume = {106},
	issn = {0003-6870},
	shorttitle = {Standards for passenger comfort in automated vehicles},
	url = {https://www.sciencedirect.com/science/article/pii/S0003687022002046},
	doi = {10.1016/j.apergo.2022.103881},
	language = {en},
	urldate = {2023-01-18},
	journal = {Applied Ergonomics},
	author = {de Winkel, Ksander N. and Irmak, Tugrul and Happee, Riender and Shyrokau, Barys},
	month = jan,
	year = {2023},
	pages = {103881},
}

@inproceedings{smit_informed_2022,
	title = {Informed sampling-based trajectory planner for automated driving in dynamic urban environments},
	doi = {10.1109/ITSC55140.2022.9922516},
	booktitle = {2022 {IEEE} 25th {International} {Conference} on {Intelligent} {Transportation} {Systems} ({ITSC})},
	author = {Smit, Robin and van der Ploeg, Chris and Teerhuis, Arjan and Silvas, Emilia},
	month = oct,
	year = {2022},
	pages = {1690--1697},
}

@inproceedings{dosovitskiy_carla_2017,
	title = {{CARLA}: {An} {Open} {Urban} {Driving} {Simulator}},
	shorttitle = {{CARLA}},
	url = {https://proceedings.mlr.press/v78/dosovitskiy17a.html},
	language = {en},
	urldate = {2023-01-16},
	booktitle = {Proceedings of the 1st {Annual} {Conference} on {Robot} {Learning}},
	publisher = {PMLR},
	author = {Dosovitskiy, Alexey and Ros, German and Codevilla, Felipe and Lopez, Antonio and Koltun, Vladlen},
	month = oct,
	year = {2017},
	pages = {1--16},
}

@article{dey_understanding_2001,
	title = {Understanding and {Using} {Context}},
	volume = {5},
	issn = {1617-4909},
	url = {https://doi.org/10.1007/s007790170019},
	doi = {10.1007/s007790170019},
	language = {en},
	number = {1},
	urldate = {2022-11-25},
	journal = {Personal and Ubiquitous Computing},
	author = {Dey, Anind K.},
	month = feb,
	year = {2001},
	pages = {4--7},
}

@misc{pribadi_grakn_2020,
	title = {The {Grakn} {Ontology}: {Simplicity} and {Maintainability}},
	shorttitle = {The {Grakn} {Ontology}},
	url = {https://blog.vaticle.com/the-grakn-ai-ontology-simplicity-and-maintainability-ab78340f5ff6},
	language = {en},
	urldate = {2022-11-25},
	journal = {Medium},
	author = {Pribadi, Haikal},
	month = mar,
	year = {2020},
}

@inproceedings{van_der_ploeg_long_2022,
	title = {Long {Horizon} {Risk}-{Averse} {Motion} {Planning}: {A} {Model}-{Predictive} {Approach}},
	shorttitle = {Long {Horizon} {Risk}-{Averse} {Motion} {Planning}},
	doi = {10.1109/ITSC55140.2022.9921750},
	booktitle = {2022 {IEEE} 25th {International} {Conference} on {Intelligent} {Transportation} {Systems} ({ITSC})},
	author = {van der Ploeg, Chris and Smit, Robin and Teerhuis, Arjan and Silvas, Emilia},
	month = oct,
	year = {2022},
	pages = {1141--1148},
}
\end{document}